\documentclass[graybox]{svmult}

\usepackage{mathptmx}       % selects Times Roman as basic font; note that this font does
\usepackage{helvet}         % selects Helvetica as sans-serif font
\usepackage{courier}        % selects Courier as typewriter font
\usepackage{graphicx}              % standard LaTeX graphics tool when including figure files

\usepackage[caption=false]{subfig} % do NOT use the subcaption package as it clashes with the Springer caption format
\usepackage{siunitx}
\usepackage{amsbsy}
\usepackage{amsmath,bm}
\usepackage{amssymb}
\usepackage{amsfonts}
\usepackage{multirow}
\usepackage[bottom]{footmisc}% places footnotes at page bottom
\usepackage{cite}
\usepackage{url}

\begin{document}

%\title*{Parallel-in-time Adjoint Sensitivity Analysisfor Time Periodic Circuits}
\title*{Periodic Adjoint Sensitivity Analysis}
% Use \titlerunning{Short Title} for an abbreviated version of
% your contribution title if the original one is too long
\titlerunning{Periodic Adjoint Sensitivity Analysis}
\author{Julian Sarpe %\inst{1}
\and
Andreas Klaedtke %\inst{2}
\and
Herbert De Gersem %\inst{3}
}
% Use \authorrunning{Short Title} for an abbreviated version of
% your contribution title if the original one is too long
\institute{Julian Sarpe, Herbert De Gersem \at Technische Universit\"at Darmstadt, Schlo{\ss}gartenstr. 8, 64289 Darmstadt, Germany
\email{julian_johannes.buschbaum@tu-darmstadt.de, herbert.degersem@tu-darmstadt.de}
  \and Andreas Klaedtke \at Robert Bosch GmbH, Robert-Bosch-Campus 1, 71272 Renningen, Germany
  \email{andreas.klaedtke@de.bosch.com}
}

%
% Use the package "url.sty" to avoid
% problems with special characters
% used in your e-mail or web address
%
\maketitle

\abstract*{%
	This paper proposes the utilization of a periodic Parareal with a periodic coarse problem to efficiently perform adjoint sensitivity analysis for the steady state of time-periodic nonlinear circuits. 
	In order to implement this method, a modified formulation for adjoint sensitivity analysis based on a transient approach is derived.
}

\abstract{%
	This paper proposes the utilization of a periodic Parareal with a periodic coarse problem to efficiently perform adjoint sensitivity analysis for the steady state of time-periodic nonlinear circuits. 
	In order to implement this method, a modified formulation for adjoint sensitivity analysis based on the transient approach is derived.
}

\section{Introduction}
Sensitivity analysis for time-periodic systems is an ongoing challenge in research until this day.
All sensitivity analysis methods require a numerical description of the system at hand.
Linear time-periodic systems can be efficiently analyzed in frequency domain.
If nonlinear devices are present, this is not possible by default.
Consequently, the analysis of nonlinear systems requires approximative methods in frequency domain, such as harmonic balance~\cite{sarpe:harmonicBalance, sarpe:TFHA} or the approximation approach used in the Plecs software~\cite{sarpe:plecs}.
\par
More often, nonlinear systems are analyzed in time-domain using transient simulations.
The nonlinearities are resolved by invoking Newton's method in each timestep.
%This can easily be performed by using a Newton iteration to find the solution for every timestep.
\par
Modified nodal analysis (MNA) is generally the method of choice for circuit analysis~\cite{sarpe:HoRuehliMNA}.
In the transient case, MNA describes the circuit in terms of a system of differential algebraic equations (DAEs) of the form
\begin{equation}\label{sarpe:eq_MNA}
	\mathbf{J}_\mathrm{C}(\dot{\mathbf{x}}(t), t) \dot{\mathbf{x}}(t) + \mathbf{J}_\mathrm{G}(\mathbf{x}(t), t) \mathbf{x}(t) = i_\mathrm{s} (t).
\end{equation}
Here, $\mathbf{J}_\mathrm{C}(\dot{\mathbf{x}}(t), t)$ gathers contributions from capacitances and inductances, $\mathbf{J}_\mathrm{G}(\mathbf{x}(t), t)$ gathers contributions from resistances, $i_\mathrm{s} (t)$ gathers the independent sources and $\mathbf{x}(t)$ contains the degrees of freedom (DoFs).
While transient simulation is the preferred approach when the entire transient must be simulated, this analysis poses a simulation overhead if only the steady state is required.
This simulation overhead is large if the transient effects (delays, relaxation) are magnitudes longer than one period.
\par
Existing methods for adjoint sensitivity analysis (ASA) in time domain require a simulation of the entire transient problem, even if only the sensitivity of the steady state must be analyzed~\cite{sarpe:nikolovaAdjoint}.
This causes simulation overhead, which the proposed method in this paper solves by modification of the transient ASA to tackle periodic problems.
\par
This paper ist structured as follows. 
First a brief description of ASA is given.
Then, Parareal and its periodic extension periodic Parareal with a periodic coarse problem (PP-PC) are presented as a possible approach to efficiently solve periodic problems.
This approach is combined with the sensitivity analysis to obtain the novel periodic ASA approach.
Lastly, a short summary and outlook are given.

\section{Adjoint Sensitivity Analysis}
\label{sarpe:sec_general}
ASA is the method of choice if the number of design parameters is larger than the number of quantities of interest (QoI).
\par
ASA was originally introduced in 1969 for transient circuit analysis~\cite{sarpe:adjointSens}. 
When analyzing time-dependent sensitivities, ASA can lag behind DSA because ASA requires each time instance to be handled as a distinct QoI. 
As a result, the choice between the two methods is a trade-off, weighing the number of parameters against the number of QoIs.
\par
To derive the ASA, the gradient of the QoIs w.r.t. the design parameters $p$ must be derived first.
It can be calculated by various methods, one of which involves employing symbolic differentiation.
Considering $i_\mathrm{s} (t)$ to be independent of $p$, the symbolic derivative of Eq.~\ref{sarpe:eq_MNA} w.r.t. $p$ reads
\begin{equation}\label{sarpe:eq_DSA}
	\mathbf{J}_\mathrm{C}(\dot{\mathbf{x}}(t), t) \frac{\mathrm{d} \dot{\mathbf{x}}(t)}{\mathrm{d} p} + \mathbf{J}_\mathrm{G}(\mathbf{x}(t), t) \frac{\mathrm{d} \mathbf{x}(t)}{\mathrm{d} p}
	\\+ \frac{\mathrm{d} \mathbf{J}_\mathrm{C}(\dot{\mathbf{x}}(t), t)}{\mathrm{d} p}\dot{\mathbf{x}}(t) + \frac{\mathrm{d} \mathbf{J}_\mathrm{G}(\mathbf{x}(t), t)}{\mathrm{d} p}\mathbf{x}(t) = 0.
\end{equation}
This formulation is usually referred to as direct sensitivity analysis (DSA)~\cite{sarpe:nikolovaAdjoint}.
ASA can be derived based on the DSA formulation~\eqref{sarpe:eq_DSA}.
\par
ASA is derived from the DSA by integrating differential equation~\eqref{sarpe:eq_DSA} over time and introducing a Lagrange multiplier $\mathbf{\lambda}$:
\begin{multline}
	\int_{t_\mathrm{0}}^{t_\mathrm{end}} \mathbf{\lambda}^\mathsf{T} (t, t_\mathrm{end})\bigg(\mathbf{J}_\mathrm{C}(\dot{\mathbf{x}}(t), t) \frac{\mathrm{d} \dot{\mathbf{x}}(t)}{\mathrm{d} p} + \mathbf{J}_\mathrm{G}(\mathbf{x}(t), t) \frac{\mathrm{d} \mathbf{x}(t)}{\mathrm{d} p}
	\\+ \frac{\mathrm{d} \mathbf{J}_\mathrm{C}(\dot{\mathbf{x}}(t), t)}{\mathrm{d} p}\dot{\mathbf{x}}(t) + \frac{\mathrm{d} \mathbf{J}_\mathrm{G}(\mathbf{x}(t), t)}{\mathrm{d} p}\mathbf{x}(t)\bigg)~\mathrm{d} t = 0.
\end{multline}
Using integration by parts, $\mathrm{d} \dot{\mathbf{x}}(t)/\mathrm{d} p$ vanishes. 
The initial and end values of $\mathrm{d} \mathbf{x}(t)/\mathrm{d} p$ must be provided explicitly:
\begin{multline}
	\int_{t_\mathrm{0}}^{t_\mathrm{end}} \left( - \dot{\mathbf{\lambda}}^\mathsf{T} (t, t_\mathrm{end}) \mathbf{J}_\mathrm{C}(\dot{\mathbf{x}}(t), t)  + \mathbf{\lambda}^\mathsf{T} (t, t_\mathrm{end})\mathbf{J}_\mathrm{G}(\mathbf{x}(t), t)  \right)\frac{\mathrm{d} \mathbf{x}(t)}{\mathrm{d} p} ~\mathrm{d} t
	\\= - \int_{t_\mathrm{0}}^{t_\mathrm{end}} \mathbf{\lambda}^\mathsf{T} (t, t_\mathrm{end})\left(\frac{\mathrm{d} \mathbf{J}_\mathrm{C}(\dot{\mathbf{x}}(t), t)}{\mathrm{d} p}\dot{\mathbf{x}}(t) + \frac{\mathrm{d} \mathbf{J}_\mathrm{G}(\mathbf{x}(t), t)}{\mathrm{d} p}\mathbf{x}(t)\right)~\mathrm{d} t 
	\\- \left[ \mathbf{\lambda}^\mathsf{T} (t, t_\mathrm{end})\mathbf{J}_\mathrm{C}(\dot{\mathbf{x}}(t), t) \frac{\mathrm{d} \mathbf{x}(t)}{\mathrm{d} p}  \right]_{t_\mathrm{0}}^{t_\mathrm{end}},
	\label{sarpe:eq_between}
\end{multline}
with the boundary conditions
\begin{subequations}
	\begin{align}\label{sarpe:eq_init}
		\frac{\mathrm{d} \mathbf{x}}{\mathrm{d} p} (t={t_\mathrm{0}}) &= 0;
		\\
		\label{sarpe:eq_terminal}
		\mathbf{\lambda} (t={t_\mathrm{end}}, t_\mathrm{end}) &= 0.
	\end{align}
\end{subequations}
The Lagrange multiplier $\mathbf{\lambda}$ is chosen to be the solution of the adjoint equation system
\begin{equation}\label{sarpe:eq_adjoint}
	\mathbf{J}_\mathrm{C}^\mathsf{T} (\dot{\mathbf{x}}(t), t) \dot{\mathbf{\lambda}} (t, t_\mathrm{end}) - \mathbf{J}_\mathrm{G}^\mathsf{T} (\mathbf{x}(t), t) \mathbf{\lambda} (t, t_\mathrm{end}) = \frac{\partial u(\mathbf{x}(t))}{\partial \mathbf{x}(t)} .
\end{equation}
Due to the boundary conditions~\eqref{sarpe:eq_init} and \eqref{sarpe:eq_terminal}, both $\mathbf{\lambda} (t, t_\mathrm{end})$ and  $\mathbf{x}(t)$ must be solved for $t \in [t_0, t_\mathrm{end}]$.
To explicitly ensure the terminal condition~\eqref{sarpe:eq_terminal} at ${t_\mathrm{end}}$, Eq.~\eqref{sarpe:eq_adjoint} must be solved backwards in time.
Both the requirement to solve for the entire time interval $t \in [t_0, t_\mathrm{end}]$ and the backwards-in-time integration are tackled by the novel approach proposed in section~\ref{sarpe:sec_periodicAdjoint}.
\par
Substituting~\eqref{sarpe:eq_adjoint} into~\eqref{sarpe:eq_between} gives
\begin{multline}
	- \int_{t_\mathrm{0}}^{t_\mathrm{end}} \left(\frac{\partial u(\mathbf{x}(t))}{\partial \mathbf{x}(t)}\right)^\mathsf{T} \frac{\mathrm{d} \mathbf{x}(t)}{\mathrm{d} p} ~\mathrm{d} t 
	\\=\int_{t_\mathrm{0}}^{t_\mathrm{end}} \frac{\mathrm{d} u}{\mathrm{d} p} (t) ~\mathrm{d} t = - \int_{t_\mathrm{0}}^{t_\mathrm{end}} \mathbf{\lambda}^\mathsf{T} (t, t_\mathrm{end})\left(\frac{\mathrm{d} \mathbf{J}_\mathrm{C}(\dot{\mathbf{x}}(t), t)}{\mathrm{d} p}\dot{\mathbf{x}}(t) + 	\frac{\mathrm{d} \mathbf{J}_\mathrm{G}(\mathbf{x}(t), t)}{\mathrm{d} p}\mathbf{x}(t)\right)~\mathrm{d} t, %- \left[ \mathbf{\lambda}^T \mathbf{J}_\mathrm{C} \frac{\mathrm{d} \mathbf{x}}{\mathrm{d} p}  \right]_{t_\mathrm{0}}^{t_\mathrm{end}}
\end{multline}
which is the integral of the sensitivity over time and the result of the transient ASA.

\section{Parareal and PP-PC}
\label{sec:parareal}
Initially introduced by Maday and Turinici in 2000~\cite{sarpe:lionsMadayParareel}, Parareal draws its foundation from multiple shooting methods~\cite{sarpe:ganderParareal}. 
One significant benefit of Parareal is its non-intrusiveness.
%, allowing it to be potentially applicable to a wide range of time-integration problems.
\par
Each iteration of Parareal, employs two distinct time domain solvers.
Firstly, a fast coarse solver (denoted as $\mathcal{G}$) gives a rough approximation of the solution $X$.
This solver is utilized for updating the initial values of $N$ temporal subintervals within each Newton iteration.
Secondly, multiple fine solvers (denoted as $\mathcal{F}$) operate in parallel over the $N$ temporal subintervals, giving a precise solution within each subinterval.
Given the non-intrusive nature of Parareal, the fine and coarse solvers can be chosen according to the problem at hand~\cite{sarpe:lionsMadayParareel, sarpe:ganderParareal}.
\par
One possible approach to extend Parareal for periodic problems is the periodic Parareal with a periodic coarse problem (PP-PC)~\cite{sarpe:ganderPPPC}.
The PP-PC differs from Parareal mainly in the definition of the update, which is extended by an additional update equation for the initial values
\begin{subequations}
	\begin{align}\label{sarpe:eq_pppcUpdate}
		\mathbf{X}^{k+1}_{n} &= \mathcal{F}(\mathbf{X}^{k}_{n-1}) + \mathcal{G}(\mathbf{X}^{k+1}_{n-1}) - \mathcal{G}(\mathbf{X}^{k}_{n-1}) \\
		\mathbf{X}^{k+1}_\mathrm{0} &= \mathcal{F}(\mathbf{X}^{k}_{N-1}) + \mathcal{G}(\mathbf{X}^{k+1}_{N-1}) - \mathcal{G}(\mathbf{X}^{k}_{N-1}).
	\end{align}
\end{subequations}
Parareal in combination with ASA was previously published~\cite{sarpe:PINTadjoint}.
In this article, this approach is extended to accommodate periodic problems as well.

\section{Periodic Adjoint Sensitivity Analysis}
\label{sarpe:sec_periodicAdjoint}
In literature, ASA for periodic systems typically involves subtracting the integrated transient ASA up to the start of the first steady state period from the integrated transient ASA up to the end of that same period~\cite{sarpe:okamotoAdjoint}:
\begin{multline}\label{sarpe:eq_periodicSensLiterature}
	\frac{\mathrm{d} U}{\mathrm{d} p} = \int_{(t_\mathrm{end}-T_\mathrm{p})}^{t_\mathrm{end}} \frac{\mathrm{d} u}{\mathrm{d} p}(t) ~\mathrm{d} t = \int_{0}^{t_\mathrm{end}} \mathbf{\lambda}^\mathsf{T}(t, t_\mathrm{end}) \left(\frac{\mathrm{d} \mathbf{J}_\mathrm{C}}{\mathrm{d} p}\dot{\mathbf{x}}(t) + \frac{\mathrm{d} \mathbf{J}_\mathrm{G}}{\mathrm{d} p}\mathbf{x}(t)\right)~\mathrm{d} t \\- \int_{0}^{(t_\mathrm{end}-T_\mathrm{p})} \mathbf{\lambda}^\mathsf{T}(t, t_\mathrm{end}-T_\mathrm{p}) \left(\frac{\mathrm{d} \mathbf{J}_\mathrm{C}}{\mathrm{d} p}\dot{\mathbf{x}}(t) + \frac{\mathrm{d} \mathbf{J}_\mathrm{G}}{\mathrm{d} p}\mathbf{x}(t)\right)~\mathrm{d} t.
\end{multline}
This procedure employs long integration intervals and requires the calculation of the circuit solution throughout the entire transient phase. 
Particularly in systems where the duration of the transient effects is magnitudes longer than the steady state period, this approach imposes substantial simulation overhead.
\par
Using the PP-PC, the circuit simulation can be significantly accelerated, reducing the simulation overhead.
To harness the benefits of this accelerated solution, adjustments must be made to accommodate periodic systems within the ASA framework.
Traditionally, the transient ASA necessitates the initial value of the circuit solution $\mathbf{x}(t=0) = \mathbf{0}$ and the terminal value of the adjoint solution $\mathbf{\lambda}(t, t_\mathrm{end}) = \mathbf{0}$.
Given that the conditions of zero initial and terminal values cannot be explicitly implemented with the PP-PC approach, certain considerations must be taken.
\par
Modifying the transient ASA for periodic problems, the integral on the right hand side is adapted to span a single period with a period duration of $T_\mathrm{p}$:
\begin{multline}\label{sarpe:eq_periodicSensIntegral}
	\frac{\mathrm{d} U}{\mathrm{d} p} = \int_{(t_\mathrm{m}-T_\mathrm{p})}^{t_\mathrm{m}} \frac{\mathrm{d} u}{\mathrm{d} p}(t) ~\mathrm{d} t = \int_{(t_\mathrm{m}-T_\mathrm{p})}^{t_\mathrm{m}} \mathbf{\lambda}^\mathsf{T}(t, t_\mathrm{end}) \left(\frac{\mathrm{d} \mathbf{J}_\mathrm{C}}{\mathrm{d} p}\dot{\mathbf{x}}(t) + \frac{\mathrm{d} \mathbf{J}_\mathrm{G}}{\mathrm{d} p}\mathbf{x}(t)\right)~\mathrm{d} t 
	\\- \left[ \mathbf{\lambda}^\mathsf{T}(t, t_\mathrm{end}) \mathbf{J}_\mathrm{C} (\mathbf{x}(t)) \frac{\mathrm{d} \mathbf{x}}{\mathrm{d} p} (t) \right]_{(t_\mathrm{m}-T_\mathrm{p})}^{t_\mathrm{m}}.
\end{multline}
In this formulation, the boundary term on the right hand side can be eliminated if
\begin{equation}\label{sarpe:eq_PppcCondition}
	\mathbf{\lambda}^\mathsf{T}(t_\mathrm{m}, t_\mathrm{end}) \mathbf{J}_\mathrm{C} (\mathbf{\dot{x}}(t_\mathrm{m})) \frac{\mathrm{d} \mathbf{x}}{\mathrm{d} p} (t_\mathrm{m}) = \mathbf{\lambda}^\mathsf{T}( t_\mathrm{m}-T_\mathrm{p}, t_\mathrm{end}) \mathbf{J}_\mathrm{C} (\mathbf{\dot{x}}(t_\mathrm{m}-T_\mathrm{p})) \frac{\mathrm{d} \mathbf{x}}{\mathrm{d} p} (t_\mathrm{m}-T_\mathrm{p}).
\end{equation}
This condition holds true only if the adjoint solution as well as the circuit solution and the system Jacobian matrices $\mathbf{J}_\mathrm{C}$ are periodic. 
Choosing $t_\mathrm{m}$ to be a time instance where both the forward and the adjoint solution are in a periodic state, the boundary terms vanish. 
With this modification, ASA can be directly utilized for periodic problems without the necessity to simulate the transient.

\section{Circuit Example}
\begin{figure}[!t]
	\centering
	\includegraphics[width=.7\linewidth]{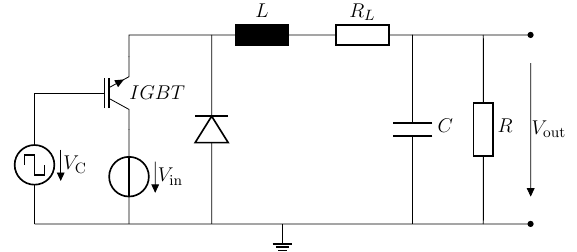}
	\caption{Schematic of the DCDC buck-converter circuit. Example taken from~\cite{sarpe:gyselinckDCDC}.}
	\label{sarpe:fig_dcdcSchematic}
\end{figure}

To demonstrate the effectiveness of the proposed method, a buck converter (depicted in Fig.~\ref{sarpe:fig_dcdcSchematic}) is used as an illustrative, academic application.
This type of buck converter is known from different publications~\cite{sarpe:gyselinckDCDC}.
The numerical simulations are executed with the model parameters collected in table~\ref{sarpe:tab_DCDCmodelParameters}.
\begin{table}[!t]
	\centering
	\caption{Model parameters used for the circuit simulation and analysis of the buck converter (Fig.~\ref{sarpe:fig_dcdcSchematic}).}
	\label{sarpe:tab_DCDCmodelParameters}
	\begin{tabular}{| l | l | l | l | l | l |} 
		\hline
		$V_\mathrm{in}$ & $f_\mathrm{s}$ & $L$ & $R_L$ & $C$ & $R$ \\
		\hline
		$\SI{100}{\volt}$ & $\SI{500}{\hertz}$ & $\SI{1}{\milli\henry}$ & $\SI{10}{\milli\ohm}$ & $\SI{100}{\micro\farad}$ & $\SI{0.8}{\Omega}$ \\
		\hline
	\end{tabular}
\end{table}
\par
\begin{figure}[!b]
	\centering
	\includegraphics[width=.95\linewidth]{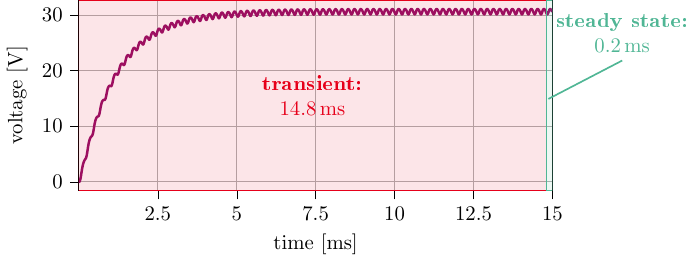}
	\caption{Output voltage $V_\mathrm{out}$ of the DCDC-Converter (Fig.~\ref{sarpe:fig_dcdcSchematic}) plotted as a function of time.}
	\label{sarpe:fig_DCDCConverterVout}
\end{figure}

The QoI is the output voltage of the buck converter.
This dynamic QoI is a charging curve with oscillations which originate from the capacitor's charging and discharging. 
Fig.~\ref{sarpe:fig_DCDCConverterVout} illustrates this QoI during the transient process as well as for the first period of steady state operation. 
The charging curve gradually reaches the desired DC value.
The calculation of the sensitivity for the first period of steady state operation using the standard approach represented by Eq.~\eqref{sarpe:eq_periodicSensLiterature} requires the circuit and the adjoint solution for the entire transient.
Considering that only the steady state information is relevant, this equates to a $\SI{98.6}{\%}$ simulation overhead.
\par
Using the novel approach introduced in Eq.~\eqref{sarpe:eq_periodicSensIntegral}, the adjoint solution is simulated using the PP-PC solver.
The solver uses the simulation parameters:
\begin{itemize}
	\item number of subintervals $N$: 2;
	\item discontinuity threshold: $10^{-4}$.
\end{itemize}
The number of subintervals had to be limited because the period is extremely small compared to the entire simulation time span.
This leads to convergence issues if the timesteps of the coarse solver are too large~\cite{sarpe:GanderConvergence}.
Alternatively, the number of subintervals could only be increased if more timesteps were used for the coarse solver, which would be detrimental to the speedup goal.
With the given simulation parameters, the PP-PC adjoint simulation converges in 9 fixed-point iterations, leading to a significant speedup even for 2 subintervals, since the timespan of the steady state is only $\SI{1.3}{\%}$ of the entire simulation time interval.
The results of the simulation are presented in table~\ref{sarpe:tab_dcdcPPPCresults}.
For reference, the Xyce direct sensitivity analysis (DSA) is used~\cite{sarpe:XyceSensitivityAnalysis}.
\par
\begin{figure}[!t]
	\centering
	\subfloat{\includegraphics[width=.4\linewidth]{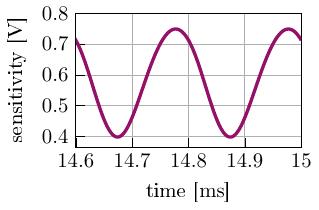}}
	\subfloat{\includegraphics[width=.4\linewidth]{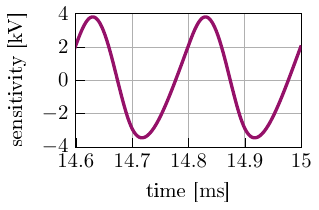}}
	\caption{Sensitivity $\mathrm{d}V_\mathrm{out}/\mathrm{d}R$ (left) and sensitivity $\mathrm{d}V_\mathrm{out}/\mathrm{d}C$ (right) of the DCDC-converter plotted as a function of time. Calculated using Xyce DSA~\cite{sarpe:XyceSensitivityAnalysis}.}
	\label{sarpe:fig_DCDCSensTDependent}
\end{figure}

\begin{table}[!t]
	\centering
	\caption{Sorted list of parameter sensitivities integrated over one period using Eq.~\eqref{sarpe:eq_periodicSensIntegral}. Periodic ASA compared with a reference solution obtained with Xyce DSA, integrated over one period. Since Xyce does not allow for sensitivity analysis w.r.t. inductances, a reference value for $L$ is not available.}
	\label{sarpe:tab_dcdcPPPCresults}
	\begin{tabular}{| c | c | c | c |}
		\hline
		device label & sensitivity value & Xyce DSA reference & error \\ \hline%\hline
		$R$ & $\SI{116.6009}{\micro\volt}$ & $\SI{116.8457}{\micro\volt}$ & $\SI{0.2095}{\%}$ \\ \hline
		$R_L$ & $\SI{-74.2263}{\micro\volt}$ & $\SI{-74.2264}{\micro\volt}$ & $\SI{0.0002}{\%}$ \\ \hline
		$L$ & $\SI{397.5767}{\nano\volt}$ & - & - \\ \hline
		$C$ & $\SI{6.358981}{\nano\volt}$ & $\SI{6.594367}{\nano\volt}$ & $\SI{3.570}{\%}$ \\ \hline
	\end{tabular}
\end{table}
The results show that the periodic approach provides accurate results for the sensitivities w.r.t. the resistances in the circuit.
Due to a high amplitude and small mean value of the oscillation of the capacitance sensitivity time series, the integral of the DSA and the periodic ASA solutions disagree.
This disagreement decreases when the DSA is resolved with more timesteps, indicating that it arises at least partially from numerical integration inaccuracies of the reference solution.
This problem is illustrated by the time dependent sensitivity in Fig.~\ref{sarpe:fig_DCDCSensTDependent}, where the mean value of the sensitivity is close to zero, but the amplitude is in the $\SI{}{\kilo\volt}$ range.
This means that depending on the design parameter, the integral of the steady state might not be sufficient to analyze the circuit behavior.

\section{Conclusions}
Periodic adjoint sensitivity analysis (ASA) provides significant computational advantages in the analysis of periodic systems compared to other sensitivity analysis approaches.
In particular, the periodic ASA is very advantageous in simulations that exhibit long transient processes, relative to the steady state period, but also contain nonlinear devices preventing an analysis in frequency domain.
Compared to transient ASA, both the circuit simulation and the subsequent sensitivity analysis are substantially sped up.
In our numerical experiments, the periodic ASA yields accurate results if the oscillation amplitude is less than its mean times $10^4$.
However, the method suffers if this is not the case, due to the discontinuity threshold convergence condition of the PP-PC, which limits the relative error to $10^{-4}$.
This leaves space for improving the accuracy of the parallel solver or the definition of the periodic integral, in order to make the method applicable to QoIs with high amplitude oscillations as well.

%%%%%%%%%%%%%%%%%%%%%%%%%%%%%%%%%%%%%%%%%%%%%%%%%%%%%%%%%%%%%%%%%%%%%%

%%%%%%%%%%%%%%%%%%%%%%%%%%%%%%%%%%%%%%%%%%%%%%%%%%%%%%%%%%%%%%%%%%%%%%

\end{document}